\documentclass[3p]{elsarticle}
\usepackage[all]{xy}
\usepackage{color}
\usepackage{url}
\usepackage{amscd, amssymb, amsfonts, amsthm, amsmath}
\usepackage{mathdots}
\usepackage{hyperref}
\usepackage{multicol}

\newcommand{\soc}{\mbox{\rm{Soc}}}
\newcommand{\topp}{\mbox{\rm{Top}}}
\newcommand{\length}{\mbox{\rm{l}}}
\newcommand{\start}{\mbox{\rm{s}}}
\newcommand{\target}{\mbox{\rm{t}}}
\newcommand{\gl}{\mbox{\rm{gldim}}}
\newcommand{\fin}{\mbox{\rm{findim}}}

\newcommand{\fidim}{\mbox{\rm{$\phi$dim}}}

\newcommand{\findim}{\mbox{\rm{findim}}}

\def\ker{\mbox{\rm{Ker}}}
\def\I{\mbox{\rm{Im}}}
\def\mod{\mbox{\rm{mod}}}
\def\ind{\mbox{\rm{ind}}}
\def\add{\mbox{\rm{add}}}

\def\pd{\mbox{\rm{pd}}}
\def\id{\mbox{\rm{id}}}

\def\ext{\mbox{\rm{Ext}}}

\def\hom{\mbox{\rm{Hom}}}
\def\enn{\hbox{\rm{End}}}

\def\rad{\hbox{\rm{rad}}}

\def\s{\mbox{\rm{s}}}
\def\t{\mbox{\rm{t}}}

\begin{document}
\newcommand{\mono}[1]{%
\gdef\puA{#1}}
\newcommand{\puA}{}
\newcommand{\faculty}[1]{%
\gdef\puC{#1}}
\newcommand{\puC}{}
\newcommand{\facultad}[1]{%
\gdef\puD{#1}}
\newcommand{\puD}{}
\newcommand{\N}{\mathbb{N}}
\newcommand{\Z}{\mathbb{Z}}
\newtheorem{teo}{Theorem}[section]
\newtheorem{prop}[teo]{Proposition}
\newtheorem{lema}[teo] {Lemma}
\newdefinition{ej}[teo]{Example}
\newtheorem{obs}[teo]{Remark}
\newtheorem{defi}[teo]{Definition}
\newtheorem{coro}[teo]{Corollary}
\newtheorem{nota}[teo]{Notation}



\title{On Lat-Igusa-Todorov algebras}

\author[add]{Marcos Barrios}
\ead{marcosb@fing.edu.uy} 

\author[add]{Gustavo Mata\corref{cor}}
\ead{gmata@fing.edu.uy}

\address[add]{Universidad de La Rep\'ublica, Facultad de Ingenier\'ia -  Av. Julio Herrera y Reissig 565, Montevideo, Uruguay}
\cortext[cor]{Corresponding Author}
\begin{abstract}
Lat-Igusa-Todorov algebras are a natural generalization of Igusa-Todorov algebras. They are defined using the generalized Igusa-Todorov functions given in \cite{BLMV} and also verify the finitistic dimension conjecture. In this article we give new ways to construct examples of Lat-Igusa-Todorov algebras. On the other hand we show an example of a family of algebras that are not Lat-Igusa-Todorov.

\end{abstract} 


\maketitle

\section{Introduction}

In an attempt to prove the finitistic dimension conjecture, Igusa and Todorov defined in \cite{IT} two functions from the objects of $\mod A$ (the category of right finitely generated modules over an Artin algebra $A$) to the natural numbers, which generalizes the notion of projective dimension. Using these functions, they showed that the finitistic dimension of Artin algebras with representation dimension at most three is finite. Nowadays, these functions are known as the Igusa-Todorov functions, $\phi$ and $\psi$.

Igusa-Todorov algebras were introduced by Wei in \cite{W1} and they verify the finitistic dimension conjecture. In particular, in the cited article, Wei proved that the class of 2-Igusa-Todorov algebras is closed under taking endomorphism algebras of projective modules. Since every Artin algebra can be realized as an endomorphism algebra of a projective and injective module over a quasi-hereditary algebra, then in case all quasi-hereditary algebra is 2-Igusa-Todorov the finitistic dimension conjecture is true.

Later, T. Conde showed, based in an article of Rouquier, that there are Artin algebras that are not Igusa-Todorov algebras (see \cite{C} and \cite{R}). 

In \cite{BLMV} Bravo, Lanzilotta, Mendoza and Vivero define the Generalized Igusa-Todorov functions and the Lat-Igusa-Todorov algebras, and prove that Lat-Igusa-Todorov algebras also verify the finitistic dimension conjecture. They also show that selfinjective algebras are Lat-Igusa-Todorov algebras, in particular the example given by T. Conde is a Lat-Igusa-Todorov algebra.

This article is organized as follows:

In Section 2, we recall the concepts given in \cite{BLMV} of $0$-Igusa-Todorov subcategories, Lat-Igusa-Todorov algebras and its properties. 

In Sections 3 and 4, we give sufficiency conditions for an algebra being a Lat-Igusa-Todorov algebra. We prove that if an algebra $A$  verifies that every module in $\Omega^n(\mod A)$ is an extension of modules of two $\mathcal{D}$-syzygy finite subcategories, then $A$ is n-Lat-Igusa-Todorov (Corollary \ref{extension D-sizigia-finita}), where $\mathcal{D}$ is a $0$-Igusa-Todorov subcategory. In particular, section 5 is dedicated to 0-Lat-Igusa-Todorov and 1-Lat-Igusa-Todorov algebras.

In Section 5, we introduce the algebras with only trivial $0$-Igusa-Todorov subcategories, i.e. every $0$-Igusa-Todorov subcategory is a subcategory of the category of projective modules. Note that: If $A$ has only trivial $0$-Igusa-Todorov subcategories, then  $A$ is an Igusa-Todorov algebra if and only if $A$ is Lat-Igusa-Todorov. We find some algebras that  have only trivial $0$-Igusa-Todorov subcategories and we also give a tool to build new family of examples (Theorem \ref{new trivial LIT}).  

Finally, Section 6 is devoted to show that some algebras are not Lat-Igusa-Todorov (Example \ref{ejemplo}).

\section{Preliminaries}

Throughout this article $A$ is an Artin algebra and $\mod A$ is the category of finitely generated right $A$-modules, $\ind A$ is the subcategory of $\mod A$ formed by all indecomposable modules, $\mathcal{P}_A \subset \mod A$ is the class of projective $A$-modules and $\mathcal{S} (A)$ is the set of isoclasses of simple $A$-modules and $A_0 = \oplus_{S \in \mathcal{S}(A)}S$. For $M\in \mod A$ we denote by $P(M)$ its projective cover, by $\Omega(M)$ its syzygy and by $\add M$ the full subcategory of $\mod A$ formed by all the sums of direct summands of $M$. For a subcategory $\mathcal{C} \subset \mod A$, we denote by $\findim (\mathcal{C})$, $\gl (\mathcal{C})$ its finitistic dimension and its global dimension respectively.

Given $A$ and $B$ algebras, if $\alpha: A \rightarrow B$ is a morphism of algebras, we know that there is an additive functor $F_{\alpha}: \mod B  \rightarrow \mod A$ such that $F_{\alpha}$ is an embedding of $\mod B$ into $\mod A$ if  $\alpha$ is an epimorphism.

A subcategory $\mathcal{C} \subset \mod A$ is add-finite if there exists $M \in \mod A$ such that $\mathcal{C} \subset \add (M)$. 

If $Q = (Q_0,Q_1,\start,\target)$ is a finite connected quiver, $\mathfrak{M}_Q$ denotes its adjacency matrix and $\Bbbk Q$ its associated path algebra. We compose paths in $Q$ from left to right. We denote by $C^n$ the oriented cycle graph with $n$ vertices. Given $\rho$ a path in $\Bbbk Q$, $\length(\rho)$, $\start(\rho)$ and $\target(\rho )$ denote the length, start and target of $\rho$ respectively.  We say that a quiver $Q$ is strongly connected if for every $v_1, v_2 \in Q_0$ there is a $\rho \in Q_1$ such that $\start(\rho) = v_1$ and $\target(\rho) = v_2$. We denote by $J$ the ideal of $\Bbbk Q$ generated by all the arrows. 


%
%

\subsection{Truncated path algebras}

We say that $A$ is a {\bf truncated path algebra} if $A = \frac{\Bbbk Q}{J^k}$ for any $k \geq 2$.
For a truncated path algebra $A$, we denote by $M^l_v(A)$ the ideal $\rho A$, where $l(\rho) = l$, $t(\rho) = v$ and $M^l(A) = \oplus_{v \in Q_0} M^l_v(A) $.

Note that if $A = \frac{\Bbbk Q}{J^k}$ is a truncated path algebra, then
$$\Omega( M_v^l(A)) = \bigoplus_{\rho:\left\{\substack{\start(\rho) = v\\ \length(\rho)= k-l}\right.} M_{ \target(\rho)}^{k-l}(A),$$
$$\Omega^2( M_v^l(A)) = \bigoplus_{\rho:\left\{\substack{\start(\rho) = v\\ \length(\rho)= k}\right.} M_{ \target(\rho)}^{l}(A).$$

For a proof of the next theorem see Theorem 5.11 of \cite{BH-ZR}, and for definitions of skeleton and $\sigma$-critical see \cite{DH-ZL}.

\begin{teo}(\cite{BH-ZR})\label{sizigia en truncadas}
Let $A$ be a truncated path algebra. If $M$ is any nonzero left $A$-module with
skeleton $\sigma$, then
$$\Omega(M) \cong \bigoplus_{\rho \text{ is } \sigma \text{-critical}} \rho A.$$
\end{teo}

Note that if $Q$ is a strongly connected quiver, then every non projective $\frac{\Bbbk Q}{J^k}$-module has infinte projective dimension.

\subsection{Igusa-Todorov functions and Igusa-Todorov algebras}
We now recall the definition of the generalized Igusa-Todorov $\phi$ function from \cite{BLMV} and some of its basic properties. 
Let us start by recalling the following version of Fitting’s Lemma.

\begin{lema}
Let $R$ be a noetherian ring. Consider a left $R$-module $M$ and $f \in \enn_R (M)$. Then, for any finitely generated $R$-submodule $X$ of $M$, there is a non-negative integer
$$\eta_{f}(X)= \min\{ k \in \mathbb{N} : f \vert _{f_m (X)} : f_m (X) \rightarrow f_{m+1} (X), \forall m \geq k\}.$$
Furthermore, for any $R$-submodule $Y$ of $X$, we have that $\eta_f(Y) \leq \eta_f(X)$.
\end{lema}

\begin{defi}
Let  $K_0(A)$ be the abelian group generated by all symbols $[M]$, with $M \in \mod A$, modulo the relations
\begin{enumerate}
  \item $[M]-[M']-[M'']$ if  $M \cong M' \oplus M''$,
  \item $[P]$ for each projective module $P$.
\end{enumerate}
\end{defi}

If $\mathcal{D} \subset \mod A$ is a subcategory such that $\mathcal{D} = \add (\mathcal{D})$ and $\Omega (\mathcal{D}) \subset \mathcal{D}$, then 
\begin{itemize}
\item The quotient group $K_{\mathcal{D}}(A) = \frac{{K}_0(A)}{\mathcal{D}}$ is a free abelian group.
\item For a subcategory $\mathcal{C} \subset \mod A$, we denote by $\langle\mathcal{C}\rangle \subset K_0(A)$ the free abelian group generated by the classes of direct summands of modules of $\mathcal{C}$ and by $[\mathcal{C}]_{\mathcal{D}} = \frac{\langle\mathcal{C}\rangle + \mathcal{D}}{\mathcal{D}}$. 
\item In particular, we denote by $\langle M \rangle = \langle \add M \rangle$ and by $\overline{ \langle M \rangle }=(\langle M \rangle + \langle \mathcal{D} \rangle )/\langle \mathcal{D} \rangle $.

\end{itemize}

\begin{lema}(\cite{BLMV})\label{lemaBLMV} Let $G$ be a free abelian group, $D$ be a subgroup of $G$, $L \in End_{\mathbb{Z}} (G)$
be such that $L(D) \subset D$ and let $k$ be a positive integer for which $L : L^k (D) \rightarrow L^{k+1}(D)$ is
a monomorphism. Then, for each finitely generated subgroup $X \subset G$, we have that

$$ \eta_L (X) \leq \eta_{\overline{L}} (\overline{X}) + k,$$

where $\overline{L} : G/D \rightarrow G/D$, $g + D  \rightarrow L(g) + D$, and $\overline{X} = (X + D)/D$.

\end{lema}

We define the {\bf Generalized Igusa-Todorov functions} as follows

\begin{defi}(\cite{BLMV})
Let $A$ be an Artin algebra and $\mathcal{D} \subset mod A$ be a subcategory such that $\Omega(\mathcal{D}) \subset \mathcal{D}$ and $\add (\mathcal{D}) = \mathcal{D}$. Let $ \bar{\Omega}_{\mathcal{D}}: K_{\mathcal{D}} \rightarrow K_{\mathcal{D}}$ be the endomorphism defined by $ \bar{\Omega}_{\mathcal{D}}([M]+\langle \mathcal{D} \rangle) = [\Omega (M)] + \langle \mathcal{D} \rangle$.  For any
$X \in mod (A)$, we set
$$\phi_{[\mathcal{D}]}(M) = \eta_{\bar{\Omega}_{\mathcal{D}}} (\overline{ \langle M \rangle }) \text{ and } \psi_{[\mathcal{D}]}(M) = \phi_{[\mathcal{D}]}(M) + \findim (\add (\Omega^{\phi_{[\mathcal{D}]}(M)}(M)))$$
where $\overline{ \langle M \rangle } =(\langle M \rangle + \langle \mathcal{D} \rangle )/\langle \mathcal{D} \rangle $. 
\end{defi}

\begin{obs}
Note that if $\mathcal{D} = \{0\}$, then $\phi_{[\mathcal{D}]} = \phi$ and $\psi_{[\mathcal{D}]} = \psi$, the Igusa-Todorov functions.
\end{obs}

Now we can define the {\bf Generalized Igusa-Todorov dimensions}.

\begin{defi}(\cite{BLMV})
Let $A$ be an Artin algebra $A$ and $\mathcal{D} \subset \mod A$ be a subcategory such that $\Omega(\mathcal{D}) \subset \mathcal{D}$ and $\add (\mathcal{D}) = \mathcal{D}$. If $\mathcal{C} \subset \mod A$, we define the {\bf $\phi_{[D]}$-dimension} and {\bf $\psi_{[D]}$-dimension} of $\mathcal{C}$, respectively, as follows:

\begin{itemize}

\item $\phi\dim_{[D]} (\mathcal{C}) = \sup \{\phi_{[D]}(M): M \in \mathcal{C}\}$,

\item $\psi\dim_{[D]} (\mathcal{C}) = \sup \{\psi_{[D]}(M): M \in \mathcal{C}\}$.

\end{itemize}

We also define the $\phi_{[\mathcal{D}]}$-dimension and $\psi_{[\mathcal{D}]}$-dimension of $A$, respectively, as follows:

\begin{itemize}

\item $\phi\dim_{[D]}(A) = \phi\dim_{[D]}(\mod A)$,

\item $\psi\dim_{[D]}(A) = \psi\dim_{[D]}(\mod A)$.

\end{itemize}

\end{defi}

The following remark summarize some propierties of the Generalized Igusa-Todorov functions.

\begin{obs}(Propositions 3.9, 3.10, and 3.12 of \cite{BLMV}) Let $A$ be an Artin algebra and $\mathcal{D} \subset \mod A$ be a subcategory such that $\Omega(\mathcal{D}) \subset \mathcal{D}$ and $\add (\mathcal{D}) = \mathcal{D}$. Then, we have the following statements, for $X, Y, M \in \mod A$.

\begin{enumerate}

\item If $M \in \mathcal{D} \cup \mathcal{P}(A)$, then $\phi_{[\mathcal{D}]} (M ) = 0$ and $\phi_{[\mathcal{D}]} (X \oplus M ) = \phi_{[\mathcal{D}]} (X)$.

\item $\phi_{[\mathcal{D}]} (X) \leq  \phi_{[\mathcal{D}]} (X \oplus Y )$ and $\psi_{[\mathcal{D}]} (X) \leq \psi_{[\mathcal{D}]} (X \oplus Y )$.

\item $\phi_{[\mathcal{D}]} \dim (\add (X))= \phi_{[\mathcal{D}]}(X)$ and $\psi_{\mathcal{[D]}} \dim(\add (X)) = \psi_{[\mathcal{D}]} (X)$.

\item $\phi_{\mathcal{[D]}} (M) \leq\phi_{\mathcal{[D]}} (\Omega(M))+1$ and $\psi_{\mathcal{[D]}} (M) \leq\psi_{\mathcal{[D]}} (\Omega(M))+1$.

\item If $Z$ is a direct summand of $\Omega^{n}(X)$ where $0 \leq t\leq\phi_{[\mathcal{D}]}(X)$ and
$\pd(Z) < \infty$, then $\pd(Z) + t \leq\psi_{[\mathcal{D}]}(X)$.

\item Suppose that $\fidim(\mathcal{D}) = 0$.

\begin{enumerate}

\item If $\pd (X) < \infty$, then $\phi_{[\mathcal{D}]} (X) = \phi(X) = \pd (X)$.

\item $\psi(X) \leq \psi_{[\mathcal{D}]} (X)$.

\item If $M \in \mathcal{D} \cup \mathcal{P}(A)$, then $\psi_{[\mathcal{D}]} (X \oplus M ) = \phi_{[\mathcal{D}]} (X)$.

\item $\psi_{[\mathcal{D}]} \dim(\mathcal{D}) = 0$.

\end{enumerate}

\end{enumerate}

\end{obs}

The following result shows the relation between the $\phi$-dimension and the $\phi_{[\mathcal{D}]}$-dimension.

\begin{teo}(\cite{BLMV})\label{BLMV} Let $A$ be an Artin algebra and $\mathcal{D} \subset \mod A$ such that $\mathcal{D} = \add (\mathcal{D})$ and $\Omega (\mathcal{D}) \subset \mathcal{D}$. Then, for every $X \in \mod A$
$$\phi(X) \leq \phi_{[\mathcal{D}]} (X) + \fidim(\mathcal{D}).$$ 
\end{teo}

\subsection{Gorenstein and Stable modules}
We denote by $ ^\bot A$ the full subcategory of $\mod A$ whose objects are those $M \in \mod A$ such that $\ext_A^i(M, A) = 0$ for $i\geq 1$.

We denote by $(\ \cdot\ )^{\ast}$ the functor $\hom_A(\ \cdot\ , A):\mod A \rightarrow \mod A^{op}$.

A finitely generated $A$-module $G$ is \textbf{Gorenstein projective} if there exists an exact sequence of $A$-modules:
$$ \xymatrix{ \ldots \ar[r] & P_{-2} \ar[r]^{p_{-2}} & P_{-1} \ar[r]^{p_{-1}} & P_0  \ar[r]^{p_0}  & P_1 \ar[r]^{p_1}& P_2 \ar[r]^{p_2} & \ldots }$$
such that $G \cong \ker (p_0) $,  $P_i$ is projective for all $i \in \mathbb{Z}$ and the following is an exact sequence:
$$ \xymatrix{ \ldots \ar[r]& {P_{2}}^{\ast} \ar[r]^{{p_{1}}^{\ast}} & {P_{1}}^{\ast} \ar[r]^{{p_{0}}^{\ast}} & {P_{0}}^{\ast} \ar[r]^{{p_{-1}}^{\ast}}  & {P_{-1}}^{\ast}\ \ar[r]^{{p_{-2}}^{\ast}}&  {P_{-2}}^{\ast} \ar[r]^{{p_{-3}}^{\ast}}& \ldots }.$$

We denote by $\mathcal{G}\mathcal{P}(A)$ the subcategory of Gorenstein projective modules.
The next properties are well known (see \cite{Zh}):

\begin{obs} Let $A$ be an Artin algebra. The following statements hold.  

\begin{enumerate}
  \item Every finite direct sum of modules of $\mathcal{G}\mathcal{P}(A)$ ($ ^\bot A$) is in $\mathcal{G}\mathcal{P}(A)$ ($ ^\bot A$)
  \item Every direct summand of modules of $\mathcal{G}\mathcal{P}(A)$ ($ ^\bot A$) is in $\mathcal{G}\mathcal{P}(A)$ ($ ^\bot A$).
  \item Every projective module is in $\mathcal{G}\mathcal{P}(A)$ ($ ^\bot A$).
  \item Every module in $\mathcal{G}\mathcal{P}(A)$ ($ ^\bot A$) is either a projective module or its projective dimension is infinite. 
\end{enumerate}
\end{obs}

Let $A$ be an algebra. We say that $A$ is a {\bf Gorenstein algebra} if $\id (A_A) <\infty$ and $\pd (D( _AA)) < \infty $.
The following results will be usefull.

\begin{prop}

Let $A$ be an Artin algebra.

\begin{enumerate}

\item If $A$ if a Gorenstein algebra, then there is a non negative integer $k$ such that  $\Omega^k (\mod A) = \mathcal{G}\mathcal{P}(A)$.

\item If $\id A_A < \infty$, then there is a non negative integer $k$ such that $\Omega^k (\mod A) = \ ^{\bot}A$.

\end{enumerate}

\end{prop}

\begin{prop} (\cite{LM}) Let $A$ be an Artin algebra, then
$$\fidim (\mathcal{G}\mathcal{P}(A)) = \fidim ( ^{\bot}A) = 0.$$

\end{prop}

\subsection{Lat-Igusa-Todorov algebras}

Lat-Igusa-Todorov algebras were introduced in \cite{BLMV} as a generalization of Igusa-Todorov algebras (see \cite{W1}). They also verify the finitistic dimension conjecture as can be seen in Theorem \ref{LIT finitista}. 

\begin{defi}
Let $A$ be an Artin algebra. If $\mathcal{D} \subset \mod A$ is a subcategory such that
\begin{enumerate}
\item $\mathcal{D} = \add (\mathcal{D})$,

\item  $\Omega (\mathcal{D}) \subset \mathcal{D}$ and

\item  $\fidim(\mathcal{D}) = 0$,
\end{enumerate}
we call it a {\bf $0$-Igusa-Todorov subcategory}.

\end{defi}

\begin{obs}\label{clases LIT}
Let $A$ be an algebra.
\begin{enumerate}
\item If $\fidim (A) = 0$, then $\mathcal{D} = \mod A$ is a $0$-Igusa-Todorov subcategory.

\item If $\fidim (A) = 1$, then $\mathcal{D} = \Omega (\mod A)$ is a $0$-Igusa-Todorov subcategory.

\item $\mathcal{G}\mathcal{P}(A)$ and $^\bot A$ are $0$-Igusa-Todorov subcategories.
\end{enumerate}
\end{obs}

\begin{defi}\label{modulo LIT}
Let $A$ be an algebra. A subcategory $ \mathcal{C} \subset \mod A$ is  ${\bf (n, V, \mathcal{D}){\text -}}${\bf Lat-Igusa-Todorov} (for short ${\bf n{\text -}LIT}$) if the following conditions are verified 
\begin{itemize}

\item There is some $0$-Igusa-Todorov subcategory $\mathcal{D} \subset \mod A$,

\item there is some $V \in \mod A$ satisfying that each $M \in \mathcal{C}$ admits an exact sequence:
$$\xymatrix{0 \ar[r]& V_1 \oplus D_1 \ar[r] & V_0 \oplus D_0 \ar[r] & \Omega^n(M)\ar[r] & 0}$$

such that $V_0, V_1 \in \add (V)$ and $D_0, D_1 \in \mathcal{D}$.
\end{itemize}
We say that $V$ is a ${\bf (n, V, \mathcal{D})\text{-}}${\bf Lat-Igusa-Todorov} {\bf module} (for short a ${\bf n\text{-}LIT}$ {\bf module}) for $\mathcal{C}$.
\end{defi}

\begin{defi}\label{algebra LIT}
We say that $A$ is a  ${\bf (n, V, \mathcal{D})\text{-}}${\bf Lat-Igusa-Todorov} {\bf algebra} (for short a  ${\bf n\text{-}LIT}$ {\bf algebra}) if $\mod A$ is $(n, V, \mathcal{D})\text{-}${\rm{LIT}}. We say that $A$ is a {\rm{LIT}} algebra if $A$ is $n\text{-}${\rm{LIT}} for some non negative integer $n$.
\end{defi}

\begin{obs}\label{Sizigia V}
Let $A$ be an algebra and $\mathcal{D}$ a $0$-Igusa-Todorov subcategory. If $V$ is a $n$-{\rm{LIT}} module, then $\Omega(V)$ is an  $(n+1)$-{\rm{LIT}} module 
\end{obs}

\begin{ej} The following are examples of {\rm{LIT}} algebras.

\begin{enumerate}

\item If $\mathcal{D} = \{0\}$, A n-LIT module $V$ is n-Igusa-Todorov and a n-LIT algebra $A$ is n-Igusa-Todorov.

\item If $\fidim (A) \leq 1$, then $A$ is a {\rm{LIT}} algebra (see Remark \ref{clases LIT}).

\item If $A$ is a Gorenstein algebra, then $A$ is a {\rm{LIT}} algebra where $\mathcal{D} = \mathcal{G}\mathcal{P}(A)$. 

\item If $\id A_A < \infty$, then then $A$ is a {\rm{LIT}} algebra where $\mathcal{D} = \ ^{\bot}A$. 

\end{enumerate}

\end{ej}

The following result show that LIT algebras verifies the finitistic dimension conjecture. For a proof see \cite{BLMV}.

\begin{teo}(\cite{BLMV})\label{LIT finitista} Let $A$ be a $(n, V, \mathcal{D})$-{\rm{LIT}} algebra. Then
$$\fin (A) \leq \psi_{[\mathcal{D}]} (V) + n + 1 < \infty.$$
\end{teo}

\section{LIT algebras and $\mathcal{D}$-syzygy finite subcategories}

In this section we show that some algebras are LIT algebras under certain properties.

\begin{obs}Let $A$ be an Artin algebra, $\mathcal{D}$ a $0$-Igusa-Todorov subcategory and $\mathcal{C} \subset \mod A$ a subcategory. If ${[\Omega^k(\mathcal{C})]}_{\mathcal{D}}$ is finitely generated, then ${[\Omega^{k+1}(\mathcal{C})]}_{\mathcal{D}}$ is finitely generated.
\end{obs}

\begin{defi}
Let $A$ an Artin algebra and $\mathcal{D}$ a $0$-Igusa-Todorov subcategory. We say that a subcategory $\mathcal{C} \subset \mod A$ is {\bf $\mathcal{D}$-syzygy finite} if ${[\Omega^k(\mathcal{C})]}_{\mathcal{D}}$ is finitely generated for some $k \in \mathbb{N}$.
\end{defi}

The following result generalizes Proposition 2.5 of \cite{W1}.

\begin{prop}\label{D-sizigia-finita}
Let $A$ be an Artin algebra and $\mathcal{D}$ be a $0$-Igusa-Todorov subcategory. If $\mod A$ is $\mathcal{D}$-syzygy finite, then $A$ is a LIT algebra.
\begin{proof}
Suppose that ${[\Omega^n(\mod A)]}_{\mathcal{D}}$ is finitely generated, then $\Omega^n (\mod A)$ is add-finite. Consider $M = \oplus M^j_i$ where $M^j_i$ is an indecomposable direct summand of $N_i \in \Omega^n(\mod A)$ for all $i, j$ that does not belong to $\mathcal{D}$. Then $M$ is a $n$-LIT module. 
\end{proof}

\end{prop}

\begin{prop}\label{sec}
Let $A$ be an Artin algebra and $\mathcal{D} \subset \mod A$ a $0$-Igusa-Todorov subcategory. If $\mathcal{C}_1$, $\mathcal{C}_2$, $\mathcal{E}$ are three subcategories of $A$-modules such that, for any $E \in \mathcal{E}$, there is an exact sequence $0 \rightarrow C_1 \rightarrow C_2 \rightarrow E \rightarrow 0$ with $C_i \in \mathcal{C}_i$ for $i = 1,2$, the next statements follows. 

\begin{enumerate}

\item If $\mathcal{C}_1$ and $\mathcal{C}_2$ are $\mathcal{D}$-syzygy finite, then $\mathcal{E}$ is $n$-{\rm{LIT}} for some $n \in \mathbb{N}$. 

\item If $\mathcal{C}_1$ is ${\mathcal{D}}$-syzygy finite and $\gl (\mathcal{C}_2) < \infty$, then $\mathcal{E}$ is $\mathcal{D}$-syzygy finite.   

\item If $\mathcal{C}_1$ is $n$-{\rm{LIT}} and  $\gl (\mathcal{C}_2) < \infty$, then $\mathcal{E}$ is $(n+1)$-{\rm{LIT}}.

\end{enumerate}

\begin{proof}

For $E \in \mathcal{E}$ there is a short exact sequence $0 \rightarrow C_1 \rightarrow C_2 \rightarrow E \rightarrow 0$ with $C_i \in \mathcal{C}_i$ for $i =1,2$. Thus, for any $n \in \mathbb{N}$ we obtain a short exact sequence $0 \rightarrow \Omega^n(C_1) \rightarrow \Omega^n(C_2) \oplus P  \rightarrow \Omega^n(E) \rightarrow 0$ for some projective $P$.

\begin{enumerate}
\item Since ${[\Omega^n(\mathcal{C}_1)]}_{\mathcal{D}}$ and ${[\Omega^n(\mathcal{C}_2)]}_{\mathcal{D}}$ are finitely generated for $n\in \mathbb{N}$, there are modules $U = \oplus_{i = 1}^t U_i$ and $V = \oplus_{j = i}^s V_j$ such that if $M_1 \in \Omega^n(\mathcal{C}_1)$ and $M_2 \in \Omega^n(\mathcal{C}_2)$, then $M_1 = \oplus_{i = 1}^t U^{\alpha_i}_i \oplus D_1$ and $M_2 = \oplus_{j = 1}^s V^{\beta_j}_j \oplus D_2$, where $D_i \in \mathcal{D}$ for $i = 1,2$ and $\alpha_i,\beta_j \in \mathbb{N}$. Hence for every $E \in \mathcal{E}$ there is a short exact sequence
$$0 \rightarrow U_1 \oplus D'_1 \rightarrow V_1 \oplus D'_2 \oplus P \rightarrow \Omega^n(E) \rightarrow 0$$
with $U_1 \in \add (U)$, $V_1 \in \add (V)$, $D_i \in \mathcal{D}$ for $i = 1,2$ and $P$ a projective module. We conclude that $\mathcal{E}$ is $n$-LIT with LIT module $U \oplus V \oplus A$.   

\item Take $n \in \mathbb{N}$ such that ${[\Omega^n(\mathcal{C}_1)]}_{\mathcal{D}}$ is finitely generated and $\gl (\mathcal{C}_2) \leq n$. Then $\Omega^n(C_2)$ is projective for every $C_2 \in \mathcal{C}_2$. It follows that $ \Omega^n(C_1) = \Omega^{n+1}(E) \oplus P$ for some projective $P$. Now we easily obtain that ${[\Omega^{n+1}(\mathcal{E})]}_{\mathcal{D}}$ is finitely generated. 

\item Take $n$ to be an integer such that $\mathcal{C}_1$ is $n$-LIT and $\gl (\mathcal{C}_2) \leq n$.
Similarly to the proof of item $(2)$, we obtain that $\Omega^n (C_1) = \Omega^{n+1} (E)  \oplus P$ for some
projective $P$. Note that there is an exact sequence $0 \rightarrow V_1 \oplus D_1 \rightarrow V_0 \oplus D_0 \rightarrow \Omega^n (C) \rightarrow 0$ with $V_i \in \add (V)$ and $D_i \in \mathcal{D}$ for $i = 0,1$, where $V$ is a $n$-LIT module. Since $P$ is projective, we can also obtain an exact sequence $0 \rightarrow V'_1 \oplus D'_1  \rightarrow V'_0 \oplus D'_0  \rightarrow \Omega^{n+1} (E) \rightarrow 0$ with $V'_i \in \add (V)$ and $D_i \in \mathcal{D}$ for $i = 0,1$. It follows that $\mathcal{E}$ is $(n+1)$-LIT with $V$ a $(n+1)$-LIT module.
\end{enumerate}
\end{proof}
\end{prop}

\begin{obs}\label{eleccion de n}
Note that in part 1 of Proposition \ref{sec},  $\min\{m: [\Omega^m(\mathcal{C}_1)] \text{ and } [\Omega^m(\mathcal{C}_2)] \text{ are finitely generated}\}$ is a possible choice of $n$.  
\end{obs}

\begin{coro}Let $A$ be an Artin algebra and $\mathcal{D} \subset \mod A$ a $0$-Igusa-Todorov subcategory.
Consider $\mathcal{C}$, $\mathcal{F}$, $\mathcal{E}$ three subcategories of $A$-modules, such that $\gl (\mathcal{F}) < \infty$ and for any $E \in
\mathcal{E}$, there is an exact sequence $$0 \rightarrow C_1 \rightarrow F_0 \rightarrow \ldots \rightarrow F_n \rightarrow E \rightarrow 0$$ with $C_1 \in \mathcal{C}$ and each $ F_i\in \mathcal{F}$. If $\mathcal{C}$ is $\mathcal{D}$-syzygy-finite ($n$-{\rm{LIT}}), then $\mathcal{E}$ is $\mathcal{D}$-syzygy finite ($n$-{\rm{LIT}}).
\begin{proof}
Denote $\mathcal{E}_0 = \mathcal{C}$, and by induction, $\mathcal{E}_{i+1} = \{ M :\exists \ 0 \rightarrow C \rightarrow F \rightarrow M \rightarrow 0 \text{ with } C \in \mathcal{E}_i \text{ and } F \in \mathcal{F}\}$. Then by hypothesis and Proposition \ref{sec}, inductively we obtain that each $E_i$ is $\mathcal{D}$-syzygy finite ($n$-LIT). Note that $\mathcal{E} \subset \mathcal{E}_{n+1}$, so $\mathcal{E}$ is also $\mathcal{D}$-syzygy finite ($n$-LIT). \end{proof}
\end{coro}
 
\begin{prop}\label{class E}
Let $A$ an Artin algebra, $\mathcal{D} \subset \mod A$ a $0$-Igusa-Todorov subcategory, and two $\mathcal{D}$-syzygy finite subcategories $\mathcal{C}_1$ and $\mathcal{C}_2$. Consider $\mathcal{E} \subset \mod A$ a subcategory such that $\forall M \in \mathcal{E}$ there exists a short exact sequence $0\rightarrow C_1\rightarrow M \rightarrow C_2 \rightarrow 0$  with $C_i \in \mathcal{C}_i$ for $i = 1,2$, then $\mathcal{E}$ is $n$-{\rm{LIT}} for some $n \in \mathbb{Z}^+$.

\begin{proof}
Suppose that for $n \in \mathbb{N}$ $[\Omega^n(\mathcal{C}_1)]_{\mathcal{D}}$ and $[\Omega^n(\mathcal{C}_2)]_{\mathcal{D}}$ are finitely generated. For any $M \in \mathcal{C}$ there are $C_i \in \mathcal{C}_i$ such that $0\rightarrow C_1\rightarrow M \rightarrow C_2 \rightarrow 0$ is a short exact sequence. Consider the following pullback diagram obtained from that short exact sequence.
$$\xymatrix{ & & 0 \ar[d] & 0 \ar[d] & \\ &   & \Omega(C_2) \ar[d] \ar[r]^{1} & \Omega(C_2) \ar[d] &  \\ 0 \ar[r] & C_1 \ar[r] \ar[d]^{1} & \ar[r] C_1 \oplus P(C_2) \ar[d] & P(C_2) \ar[r] \ar[d]& 0 \\ 0 \ar[r] & C_1 \ar[r] & M \ar[d]\ar[r] & C_2 \ar[r] \ar[d] & 0 \\ & & 0 & 0 & }$$
It is easy to check that $\Omega^n(\mathcal{E})$ is $n$-LIT, just apply part 1 of Proposition \ref{sec} to the middle column in the above diagram. 
\end{proof}

\end{prop}

The following result follows directly from the previous proposition.

\begin{coro}\label{extension D-sizigia-finita}
Let $A$ an Artin algebra, $\mathcal{D}$ a $0$-Igusa-Todorov subcategory for $\mod A$. If there are two $\mathcal{D}$-syzygy finite subcategories $\mathcal{C}_1$ and $\mathcal{C}_2$ such that for every $M \in \mod A$ there is a short exact sequence 
$$0 \rightarrow C_1 \rightarrow \Omega^n (M) \rightarrow C_2\rightarrow 0$$
with $C_i \in \mathcal{C}_i$, then $A$ is a $n$-{\rm LIT} algebra.
\end{coro}

\section{Small LIT algebras}

Throughout this section, we identify $0$-LIT and $1$-LIT algebras under conditions in the category of modules, in quotients, and its categories of modules.

The first result is a generalization of Proposition 3.2 from \cite{W1}. This result allows us to identify $0$-LIT algebras.

\begin{prop}
Let $A$ be an Artin algebra and $\mathcal{D} \subset \mod A$ a $0$-Igusa-Todorov subcategory. Consider two ideals $I, J$ with $ J I = 0$. Then $A$ is a $0$-{\rm{LIT}} algebra provided that

\begin{enumerate}

\item $\ind \frac{A}{I} \setminus \mathcal{D} \subset \mod A$ and $\ind \frac{A}{J}\setminus \mathcal{D} \subset \mod A$ are finite sets.  

\item $\ind \frac{A}{I} \setminus \mathcal{D} \subset \mod A$ is finite, $\frac{A}{J}$ is projective in $\mod A$ and ${[\Omega(\mod \frac{A}{J})]}_{\mathcal{D}}$ is finitely generated.

\end{enumerate}

\begin{proof}
For any $N \in \mod A$, we have a short exact sequence $0 \rightarrow NJ \rightarrow N \rightarrow  \frac{N}{NJ} \rightarrow 0$. Note that $ (NJ)I = 0$ and $ (\frac{N}{NJ})J = 0$, so $NJ$ is also in $\mod \frac{A}{I}$ and $\frac{N}{NJ}$ is also in $\mod \frac{A}{J}$.

Consider the following pullback diagram obtained from the above short exact sequence.
$$\xymatrix{ & & 0 \ar[d] & 0 \ar[d] & \\ &   & \Omega(\frac{N}{NJ} ) \ar[d] \ar[r]^{1} & \Omega(\frac{N}{NJ}) \ar[d] &  \\ 0 \ar[r] & NJ \ar[r] \ar[d]^{1} & \ar[r] NJ \oplus P(\frac{N}{NJ}) \ar[d] & P(\frac{N}{NJ}) \ar[r] \ar[d]& 0 \\ 0 \ar[r] & NJ \ar[r] & N \ar[d]\ar[r] & \frac{N}{NJ} \ar[r] \ar[d] & 0 \\ & & 0 & 0 & }$$

Both items follow by Remark $\ref{eleccion de n}$ applied to the middle row in the diagram. 

%
%
%
\end{proof}
\end{prop}

The following two results are generalizations of Theorem 3.4 and Corollary 3.5 of \cite{W1} respectively.

\begin{prop}\label{1-LIT}
Let $A$ be an Artin algebra, $\mathcal{D} \subset \mod A$ a $0$-Igusa-Todorov subcategory and $I$ an ideal with $\rad(A)I = 0$. If $\mod \frac{A}{I} \subset \mod A$ is $0$-{\rm{LIT}}, then $A$ is a $1$-{\rm{LIT}} algebra.

\begin{proof} By hypothesis, for any $M \in \mod A$, we have that $\Omega(M)I \subset  \rad(P(M))I = 0$. Then $\Omega(M)$ is also an $\frac{A}{I}$-module. Since $\mod \frac{A}{I} \subset \mod A$ is $0$-LIT with a LIT-module $V$, then we obtain an exact sequence of $A$-modules $0 \rightarrow V_1 \oplus D_1 \rightarrow V_0 \oplus D_0 \rightarrow \Omega(M) \rightarrow 0$ with $V_0 , V_1 \in \add (V)$ and $D_0, D_1 \in \mathcal{D}$. Hence, we conclude that $A$ is a $1$-LIT algebra with a LIT module $V$.
\end{proof}

\end{prop}

\begin{coro}
Let $A$ be an Artin algebra and $\mathcal{D} \subset \mod A$ a $0$-Igusa-Todorov subcategory. If $\rad^{2n+1}(A) = 0$ and $\ind \frac{A}{\rad^n (A)} \setminus \mathcal{D} \subset \mod A$ is finite, then $A$ is $1$-{\rm{LIT}}.
\begin{proof}
We have the following embeddings of module categories
$$\mod \frac{A}{\rad^n (A)} \subset \mod \frac{A}{\rad^{2n}(A)} \subset \mod A $$  
Consider $I = J = \frac{\rad^n A}{\rad^{2n} (A)}$ ideal of $\frac{A}{\rad^{2n}(A)}$. Observe that $IJ = 0$. If $M\in \mod \frac{A}{\rad^{2n}(A)}$, then $JM \in \mod \frac{A}{\rad^n (A)}$ and $\frac{M}{JM} \in \mod  \frac{A}{\rad^n (A)}$ and by Proposition \ref{class E} we conclude that $\mod \frac{A}{\rad^{2n}(A)} \subset \mod A$ is $0$-LIT. Finally, by Proposition \ref{1-LIT} $A$ is $1$-LIT.
\end{proof}
\end{coro}

\section{Algebras with only trivial $0$-Igusa-Todorov subcategories}

In this section we build algebras with only trivial $0$-Igusa-Todorov subcategories. We will use these results in section $7$ to construct examples of non LIT algebras.

\begin{defi}
Let $A$ be an Artin algebra. We say that $A$ has {\bf only trivial $0$-Igusa-Todorov subcategories} if for all $0$-Igusa-Todorov subcategory $\mathcal{D}$, $\mathcal{D} \subset \mathcal{P}_A$.
\end{defi}

\begin{defi}
Let $A$ be an Artin algebra. For $M \in \mod A$ we define
$$\gamma (M) = \fidim(\add\{N\text{: } N \text{ is a direct summand of } \Omega^n(M) \text{ for some }n\in \mathbb{N}\}).$$
\end{defi}

\begin{prop}
Let $A$ be an Artin algebra. The following statements are equivalent

\begin{enumerate}

\item $A$ has only trivial $0$-Igusa-Todorov subcategories.

\item $\min\{\gamma(M) \text{: such that } M \in \mod A \setminus \mathcal{P}_A\} \geq 1$.

\item  $\min\{\gamma(M) \text{: such that } M \in \ind A \setminus \mathcal{P}_A\} \geq 1$.

\end{enumerate}

\begin{proof} We prove the equivalences.

\begin{itemize}
\item $(1 \Rightarrow 2)$ Let $M \in \mod A \setminus \mathcal{P}_A$. It is clear that $\mathcal{C}_M = \{N\text{: } N \vert \Omega^n(M) \mbox{ for some } n\in \mathbb{N} \}$ verifies the first two axioms for a $0$-Igusa-Todorov subcategory. Since $A$ has only trivial $0$-Igusa-Todorov subcategories, $\fidim (\mathcal{C}_M) = \gamma(M) \geq 1$.

\item $(2 \Rightarrow 3)$ It is a particular case.

\item $(3 \Rightarrow 1)$ Let $\mathcal{D}$ be a non trivial subgategory such that is closed by syzygies and direct summands. Then there is a non projective indecomposable module $M \in \mathcal{D}$. By hypothesis $\gamma(M) \geq 1$ so there is $N \in \mathcal{D}$ such that $\phi(N) \geq 1$. We deduce that $\mathcal{D}$ is not a $0$-Igusa-Todorov subcategory.
\end{itemize}
\end{proof}
\end{prop}

\begin{prop} \label{truncadas} The following algebras have only trivial $0$-Igusa-Todorov subcategories 

\begin{enumerate}

\item If $A = \frac{\Bbbk Q}{J^2}$ is a non selfinjective radical square zero algebra such that $Q$ is strongly connected and the adjacence matrix $\mathfrak{M}_Q$ of $Q$ is not invertible.

\item If $A = \frac{\Bbbk Q}{J^k}$ is a truncated path algebra such that $Q$ is strongly connected algebra with at least one loop and the adjacence matrix $\mathfrak{M}_Q$  of $Q$ is not invertible. 

\end{enumerate}

\begin{proof}

\begin{enumerate}

\item By Proposition 4.14 and Theorem 4.32 of \cite{LMM}, $\phi(A_0) \geq 1$. If $M \in \ind A \setminus \mathcal{P}_A$, then $\Omega(M) \subset \add (A_0)$. Since $Q$ is strongly connected quiver, $A_0$ has no projective summands. On the other hand, since $Q$ is strongly connected, then $A_0 \in \add (\oplus_{k=1}^n\Omega^k(M))$, and it follows the thesis.

\item By Remark 11 of \cite{BMR}, $\phi(M^l(A) \oplus M^{k-l}(A)) \geq 1$ for every $1\leq l \leq k-2$. If $M$ is not a projective module, then $\Omega(M) = M^l_v(A) \oplus N$ for some $1\leq l \leq k-2$, $v \in Q_0$. On the other hand, since $Q$ is strongly connected and has a loop, then $M^l(A) \oplus M^{k-l}(A) \in  \add (\oplus_{k=1}^n\Omega^k(M))$, and it follows the thesis.

\end{enumerate}
 
\end{proof} 
  
\end{prop}

The following example shows that it is necessary to have at least one loop in the case of truncated path algebras of the above proposition.

\begin{ej}

Consider the algebra $A = \frac{\Bbbk Q}{J^8}$, with $Q$ the following quiver

$$
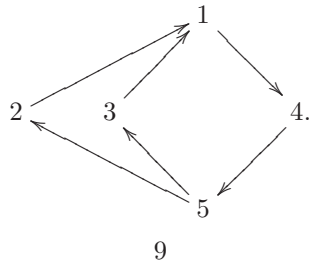
\xymatrix{ && 1 \ar[rd] & \\ 2 \ar[urr]& 3 \ar[ur]& & 4. \ar[ld]\\ && 5 \ar[ul] \ar[ull] & }$$

Let $M$ be the $A$-module given by the representation below

$$
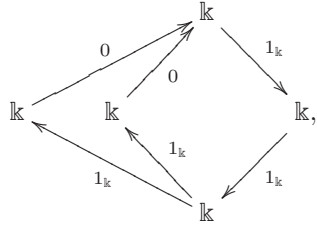
\xymatrix{ && \Bbbk \ar[rd]^{1_{\Bbbk}} & \\ \Bbbk \ar[urr]^{0} & \Bbbk \ar[ur]_0& & \Bbbk, \ar[ld]^{1_{\Bbbk}}\\ && \Bbbk \ar[ul]_{1_{\Bbbk}} \ar[ull]^{1_{\Bbbk}} & }$$
then $\Omega(M) = M^2$, and $\gamma(M) = \phi (M) = 0$. We conclude that $A$ does not have only trivial $0$-Igusa-Todorov subcategories.

\end{ej}

\begin{defi}

Let $A = \frac{\Bbbk Q}{I}$ a finite dimensional algebra. If $\bar{Q}$ is a full subquiver of $Q$ and $B = \frac{\Bbbk \bar{Q}}{I\cap \Bbbk Q}$, then we denote by $\pi_B : \mod A \rightarrow \mod B$ the restriction functor.    

\end{defi}

\begin{teo}\label{new trivial LIT}

Let $A=\frac{\Bbbk Q}{I}$ a finite dimensional algebra such that there are two disjoint full subquivers $\Gamma$ and $\bar{\Gamma}$ of $Q$ which verifies:

\begin{itemize}

\item $\bar{\Gamma}$ has no sinks.

\item $Q_0 = \Gamma_0 \cup \bar{\Gamma}_0$.

\item For all $v \in \Gamma_0$ there is an arrow $\alpha_v \in Q_1$ such that $\s(\alpha_v) = v$ and $\t(\alpha_v) = w \in \bar{\Gamma}_0$. 

\item There are no arrows $\alpha \in Q_1$ with $\s(\alpha) \in \bar{\Gamma}_0$ and $\t(\alpha) \in \Gamma_0$.

\item For all $\alpha \in Q_1$ such that $\s(\alpha) \in \Gamma_0$ and $\t(\alpha) \in \bar{\Gamma}_0$ then $\alpha \beta = 0 = \delta \alpha$ for all $\beta, \delta \in Q_1$.

\end{itemize}

If $B = \frac{\Bbbk \bar{\Gamma}}{I \cap \Bbbk \bar{\Gamma}}$ has only trivial $0$-Igusa-Todorov subcategories, then $A$ has only trivial $0$-Igusa-Todorov subcategories.

\begin{proof}

Let $B$ and $C$  be the algebras $B = \frac{\Bbbk \bar{\Gamma}}{I \cap \Bbbk \bar{\Gamma}}$ and $C = \frac{\Bbbk \Gamma}{I \cap \Gamma}$ respectively. It is easy to see that $\Omega (\mod A) \subset \mod B \oplus \mod C \oplus \{\oplus P_v : v \in \Gamma_0\}$. 
Notice that $\mod B$ has no simple projective modules.
Consider $\mathcal{D}$ a $0$-Igusa-Todorov subcategory for $A$.\\

{\bf \underline{Claim:}} $\mathcal{D} \cap \mod B$  is a $0$-Igusa-Todorov subcategory for $B$.\\

Since $\mathcal{P}_B \subset \mathcal{P}_A$, then $\Omega_B(M) = \Omega_A(M)$ for all $M \in \mod B$. Hence $\Omega_B(M) \in \mathcal{D} \cap \mod B$ and $\phi_B(M) = \phi_A(M) = 0$ for all $M \in \mathcal{D} \cap \mod B$. On the other hand consider $M \in \mod B$, if $N$ is a direct summand of $M$ in $\mod A$, it is clear that $N \in \mod B$.

As a consequence of the claim, it is clear that for $M \in \mathcal{D}\setminus \mathcal{P}_A$, if $N$ is a direct summand of $\Omega(M)$ from $\mod B$ then $N \in \mathcal{P}_B$.

Suppose $M \in \mathcal{D}\setminus \mathcal{P}_A$, then $\Omega(M)$ is not projective. Hence $\Omega(M)$ has a direct summand in $\mod C$. Since there is a simple $B$-module $S$ such that $S$ is a direct summand of $\Omega^2(M)$, then $\Omega^2(M)$ has a non projective direct summand in $\mod B$. Finally if we apply the claim to $\Omega(M)$ is a projective module, and this is absurd. 

\end{proof}

\end{teo}

\section{Examples of non LIT algebras}

In this section, we give an example of a family of finite dimensional algebras that are not LIT.

\begin{ej}\label{ejemplo}
Let $B =\frac{\Bbbk Q}{I_B}$ be a finite dimensional $\Bbbk$-algebra and $C = \frac{\Bbbk Q'}{J^2}$, where $Q'$ is the following quiver
$$\xymatrix{ 1 \ar@(ul, dl)_{\bar{\beta}_1} \ar@/^2mm/[r]^{\beta_1} & 2 \ar@(ur, dr)^{\bar{\beta}_2} \ar@/^2mm/[l]^{\beta_2} }$$
Consider $A = \frac{\Bbbk \Gamma}{I_A}$, with

\begin{itemize}

\item $\Gamma_0 = Q_0 \cup Q'_0$,

\item $\Gamma_1 = Q_1 \cup Q'_1 \cup \{ \alpha_i : i \rightarrow 1 \text{ } \forall i \in Q_0\}$ and

\item $I_A = \langle I_B, J^2_{C} , \{ \lambda \alpha_i, \alpha_i\lambda \text{ } \forall \lambda \text{ such that }\length(\lambda)\geq 1 \} \rangle$. 

\end{itemize}

Note that 

\begin{itemize}

\item If $M \in \mod A$, then $\pd M = \left\{\begin{array}{cc} 0,& \text{or} \\ \infty. \end{array} \right.$

\item $K_1(A) \subset \langle [M]: M\in \mod B \subset \mod A \rangle \times [S_1] \times [S_2]$.

\item If $M \in \mod B$, then $\Omega_A(M) = \Omega_B(M) \oplus S_1^{\dim\topp(M)}$.

\item If $M \in \mod A$ and $\pd(M) = \infty$, then $S_1$ and $S_2$ are direct summands of $\Omega^3_A(M)$. As a consequence $A$ is a $LIT$ algebra if and only if $A$ is an Igusa-Todorov algebra (Use Theorem \ref{new trivial LIT} and Proposition \ref{truncadas}).


\end{itemize}

\end{ej}

\begin{obs}\label{remarquita}
Let $A$ be an algebra as in Example \ref{ejemplo} where $B$ is a selfinjective algebra. If $0\rightarrow V_B \oplus S \rightarrow P \rightarrow W_B \oplus \bar{S} \rightarrow 0$ is a short exact sequence in $\mod A$ with $V_B, W_B \in \mod B \setminus \mathcal{P}_B$, $P \in \mathcal{P}_A$ and $S, \bar{S} \in \add(S_1 \oplus S_2)$, then there is a short exact sequence $0\rightarrow V_B \rightarrow \bar{P} \rightarrow W_B  \rightarrow 0 $ in $\mod A$ with $\bar{P} \in \mathcal{P}_B$.
\end{obs}

\begin{obs}
Let $A$ be an algebra as in Example \ref{ejemplo} where $B$ is a selfinjective algebra. If $A$ is an 1-Igusa Todorov algebra, then $B$ is also an 1-Igusa Todorov algebra.  
\end{obs}

\begin{lema}\label{K_1(B) A autoinyectiva}
Let $A$ be an algebra as in Example \ref{ejemplo} where $B$ is a selfinjective algebra, then 
$$K_1(A) = \langle [M]: M\in \mod B \setminus \mathcal{P}_B \subset \mod A \rangle \times [S_1] \times [S_2].$$

\begin{proof}
It easy to see that $S_1, S_2 \in K_1(A)$, and If $P \in \mathcal{P}_B$ then $P \not \in K_1(A)$. On the other hand consider $V_B \in \mod B \setminus \mathcal{P}_B$. Since $B$ is a selfinjective algebra, there is a short exact sequence in $ \mod B$ as follows 
$$0 \rightarrow V_B \rightarrow P \rightarrow W_B \rightarrow 0,$$
where $P \in \mathcal{P}_B$. From the previous short exact sequence, we can construct the following short exact sequence in $\mod A$.
$$0 \rightarrow V_B \oplus S_1^{\topp (W_B)} \rightarrow \bar{P} \rightarrow W_B \rightarrow 0,$$
where $\bar{P} \in \mathcal{P}_A$. We deduce that $V_B \in K_1(A)$.
\end{proof}

\end{lema} 

As a consequence of the proof of Lemma \ref{K_1(B) A autoinyectiva} we have the next result.

\begin{coro} \label{Sizigia de Omega}
Let $A$  be an algebra as in Example \ref{ejemplo} where $B$ is a selfinjective algebra. Then the next statements follows 
\begin{enumerate}
\item If $V \in \Omega_A(\mod A)$, there is a semisimple $S \in \mod A$ such that $0 \rightarrow V \oplus S \rightarrow P \rightarrow W \rightarrow 0$, with $P \in \mathcal{P}_A$ and $W \in \Omega_A(\mod A)$.

\item $\bar{\Omega}\vert_{\langle [M]: M\in \mod B \setminus \mathcal{P}_B \rangle}$ is injective.

\end{enumerate}

\begin{proof}\

\begin{enumerate}

\item The $A$-module $V$ can be decomposed into $V = V_B \oplus S_1^{m_1} \oplus S_2^{m_2}$ with $V_B \in \mod B$.

Let $W_B$ be a preimage of $V_B$, and $\bar{W}_B$ a preimage of $W_B$ as in Lemma \ref{K_1(B) A autoinyectiva}. It is easy to see that $\Omega(S_1) = \Omega(S_2) = S_1\oplus S_2$, then 
$$\Omega (W_B \oplus S_1^{\topp (\bar{W}_B)+m_1+m_2}) = V_B \oplus  {S_1}^{ \topp (W_B) + \topp (\bar{W}_B)+ m_1+m_2} \oplus {S_2}^{\topp (\bar{W}_B)+m_1+m_2} $$ 
\item Is a direct consequence of Lemma \ref{K_1(B) A autoinyectiva}

\end{enumerate}

\end{proof}

\end{coro}

%
%

\begin{prop}
Let $A$ as in Example \ref{ejemplo} where $B$ is a selfinjective algebra. If $A$ is $m$-Igusa Todorov, then $A$ is $1$-Igusa Todorov. 

\begin{proof}
If $A$ is a $m$-Igusa Todorov algebra with $m>1$, we can assume, by Remark \ref{Sizigia V}, that there exist an Igusa Todorov module $V$ such that $V \subset \Omega_A(\mod A)$. Assume that $A_0$ is a direct summand of $V$.
Given the short exact sequences 
$$0 \rightarrow V_1 \stackrel{u_m}{\rightarrow} V_0 \stackrel{v_m}{\rightarrow} \Omega^m(M) \rightarrow 0 \text{, and } 0 \rightarrow \Omega^{m}(M) \stackrel{i_{m-1}}{\rightarrow}  P_{m-1} \stackrel{p_{m-1}}{\rightarrow}  \Omega^{m-1}(M) \rightarrow 0,$$
we can construct the following commutative diagram with exact columns and rows 
$$\xymatrix{  & 0 \ar[d] & 0 \ar[d] & 0 \ar[d] & \\  0\ar[r] &  V_1 \ar[d]^{\gamma_{m-1} u_{m}} \ar[r]^{u_m}  & V_0 \ar[d]^{\mu_{m-1}} \ar[r]^{v_m} & \Omega^{m}(M) \ar[d]^{i_{m-1}} \ar[r] & 0 \\ 0 \ar[r] & Q_{m-1}\ar[r]^{\iota_{m-1}} \ar[d] &  Q_{m-1}  \oplus P_{m-1} \ar[r]^{\Pi_{{m-1}}} \ar[d] & P_{m-1} \ar[r] \ar[d]^{p_{m-1}}& 0 \\ 0 \ar[r] & W_1  \ar[d] \ar[r]& W_0 \ar[d]\ar[r] & \Omega^{m-1}(M) \ar[r] \ar[d] & 0 \\ & 0 & 0 & 0 & }$$
where the maps $\iota_{m-1}$ and $\Pi_{m-1}$ are the canonical inclusion and projection respectively, and $\mu_{m-1} = \left(\begin{array}{c}\gamma_{m-1} \\ i_{m-1}v_{m} \end{array} \right) $. Consider $S \in \mod A$ a semisimple module and $\lambda_{m-1}: S \rightarrow Q_{m-1}$ such that $\delta_{m-1}: V_0 \oplus S \rightarrow Q_{m-1} \oplus P_{m-1} $, given by $\left(\begin{array}{cc}\gamma_{m-1} & \lambda_{m-1} \\ i_{m-1}v_{m} & 0 \end{array} \right) $, is a monomorphism and $(Soc (Q_{m-1}), 0) \subset \I (\delta_{m-1})\vert_{Soc (V_0) \oplus S}$. Consider $\epsilon_{m-1}: V_1\oplus S \rightarrow Q_{m-1}$ given by $\epsilon_{m-1} = (\gamma_{m-1}u_m, \lambda_{m-1}) $.\\

{\bf \underline{Claim:}} The map $\epsilon_{m-1}$ is a monomorphism.\\

Suppose there exist $v \in V_1$ and $s\in S$ such that $\epsilon_{m-1}(v,s) = \gamma_{m-1}u_m(v) + \lambda_{m-1}(s) = 0$. Since $u_m(v) \in V_0$ and $s\in S$, then 
$$\delta_{m-1}(u_m(v), s) = \left(\begin{array}{cc}\gamma_{m-1} & \lambda_{m-1} \\ i_{m-1}v_{m} & 0 \end{array} \right) \left(\begin{array}{c} u_m(v) \\ s \end{array}\right) = (\gamma_{m-1}u_m(v) + \lambda_{m-1}(s),  i_{m-1}v_{m} u_m(v) ) = (0,0)$$
Since $\delta_{m-1}$ and $u_m$ are monomorphisms, then $v = 0$ and $s= 0$.\\  
  
From the above diagram and the maps $\epsilon_{m-1}$, $\lambda_{m-1}$ and $\bar{v}_m = (v_{m}, 0)$, by Lemma $3\times 3$, we obtain the following diagram.  
$$\xymatrix{  & 0 \ar[d] & 0 \ar[d] & 0 \ar[d] & \\  0\ar[r] &  V_1 \oplus S \ar[d]^{\epsilon_{m-1}} \ar[r]^{u_m\oplus 1_S} & V_0 \oplus S\ar[d]^{\delta_{m-1}} \ar[r]^{\bar{v}_m} & \Omega^{m}(M) \ar[d]^{i_{m-1}} \ar[r] & 0 \\ 0 \ar[r] & Q_{m-1}\ar[r] \ar[d] &  Q_{m-1}  \oplus P_{m-1} \ar[r]^{\Pi_{{m-1}}} \ar[d]^{q_{m-1}} & P_{m-1} \ar[r] \ar[d]^{p_{m-1}}& 0 \\ 0 \ar[r] & \bar{W}_1  \ar[d] \ar[r]& \bar{W}_0 \ar[d]\ar[r]^{\omega_{m-1}} & \Omega^{m-1}(M) \ar[r] \ar[d] & 0 \\ & 0 & 0 & 0 & }$$ 
We denote by $\bar{W}_0 = ((W_0)^i, T_{\alpha})$,  $Q_{m-1} = ((Q_{m-1})^i, \bar{T}_{\alpha})$ and $P_{m-1} = ((P_{m-1})^i, \tilde{T}_{\alpha})$ as representations.\\ 

{\bf \underline{Claim:}} $[\bar{W}_0] \in K_1(A)$.\\

Let $w \in \bar{W}_0$ such that $w \not= 0$ and $e_1 w = w$ (the case $e_2w = w$ is easier and left to the reader). We want to prove that $w \not \in \I \sum_{\alpha: j \rightarrow 1} T_{\alpha}$ and $T_{\beta_1}(w) = T_{\bar{\beta}_1}(w) = 0$.

Suppose there exists $w' \in W_0$ such that $\sum_{\alpha: j \rightarrow 1} T_{\alpha} (w') = w$, then $\omega_{m-1}(w) = 0$. Since $q_{m-1}$ is an epimorphism, there exist $x, x' \in Q_{m-1}\oplus P_{m-1}$ where $q_{m-1}(x) = w$, $q_{m-1}(x') = w'$ and $\sum_{\alpha: j \rightarrow 1} \bar{T}_{\alpha}+\tilde{T}_{\alpha}(x') = x$. We deduce that $x \in S_1 \subset \soc(Q_{m-1}\oplus P_{m-1})$.
 
Now consider $y, y' \in P_{m-1}$ such that $\Pi_{m-1}(x) = y$ and $\Pi_{m-1}(x') = y'$, since $(Soc (Q_{m-1}), 0) \subset \I (\delta_{m-1})\vert_{Soc (V_0) \oplus S}$ it is clear that $y \not = 0$. 
By the previous diagram there is an element $z \in S_1 \subset \soc( \Omega^{m}(M))$ such that $i_{m-1}(z) = y$. 

Since $\bar{v}_m$ is an epimorphism there is an element $v \in S_1 \subset \soc (V_0)$ such that $\bar{v}_m(v) = z$. Again, by the previous diagram $\Pi_{m-1}(x-\delta_{m-1}(v) )=0$, then $x-\delta_{m-1}(v) \in Q_{m-1}$. 
Since $x, \delta_{m-1}(v) \in \soc (Q_{m-1}\oplus P_{m-1})$, it is clear that $x-\delta_{m-1}(v) \in (\soc (Q_{m-1}),0)$. Therefore there exists $v' \in \soc (V_0 \oplus S)$ such that $\delta_{m-1}(v') = x-\delta_{m-1}(v)$. 
It is an absurd since $0 = q_{m-1}\delta_{m-1}(v') = q_{m-1} (x-\delta_{m-1}(v)) = q_{m-1}(x) = w \not = 0$.
  
Now, if we suppose that $T_{\beta_1} (w) \not = 0$ ($T_{\beta_2} (w) \not = 0$). Consider $x = T_{\beta_2}(w)$ ($x = T_{\beta_1} (w)$) and the proof follows as above.

Finally, by Remark \ref{Sizigia de Omega}, there is a semisimple module $\bar{S}$ such that $\bar{W}_0 \oplus \bar{S} \in \Omega_A(\mod A)$.

From the below short exact sequence of the previous commutative diagram we build the following short exact sequence
$$0 \rightarrow \bar{W}_1 \oplus \bar{S} \rightarrow \bar{W}_0 \oplus \bar{S} \rightarrow \Omega^m(M)\rightarrow  0.$$
Since $\Omega^{m-1}(M)$ and $\bar{W}_0 \oplus \bar{S}$ belong to $\Omega_A(\mod A)$, then $\bar{W}_1 \oplus \bar{S} \in \Omega_A(\mod A)$. By Remark \ref{remarquita}, there exist $W \in \Omega_A (\mod A)$ such that $\bar{W}_1$, $\bar{W}_1$ belong to $\add(W)$ for all $M \in mod A$ and the thesis follows.\end{proof}
\end{prop}

We finally give an example of an Artin algebra that is not Lat-Igusa-Todorov.

\begin{ej}
Let $A$ as in Example \ref{ejemplo} where $B$ is a selfinjective algebra. 
If $B$ is not an Igusa-Todorov algebra (see 4.2.10 of \cite{C} and Corollary 4.4 of \cite{R}), then $A$ is not an {\rm{LIT}} algebra. However, by Theorem 5.2 of \cite{BM} $\fidim(A) \leq 3$ (in fact $\fidim (A) = 2$), and $A$ verifies the finitistic dimension conjecture.
\end{ej}

{\bf Acknowledgements.} The authors thank Professor Marcelo Lanzilotta for helpful comments and
recommendations which helped to improve the quality of the article.

\bibliographystyle{elsarticle-harv} 


\end{document}